\documentclass[11pt,oneside]{amsart}
\usepackage[naturalnames]{hyperref}
\usepackage{a4wide}
\usepackage{array}
\usepackage{graphics}
\usepackage{longtable}
\usepackage{multirow}
\usepackage[ruled,norelsize]{algorithm2e}

\theoremstyle{plain}

\newtheorem{theorem}{Theorem}[section]
\newtheorem{lemma}[theorem]{Lemma}
\newtheorem{proposition}[theorem]{Proposition}

\theoremstyle{definition}
\newtheorem{definition}[theorem]{Definition}
\theoremstyle{remark}

\newcommand{\bbr}{{\mathbb R}}
\newcommand{\bbp}{{\mathbb P}}
\newcommand{\bbc}{{\mathbb C}}  
\newcommand{\bbh}{{\mathbb H}}
\newcommand{\bbn}{{\mathbb N}}
\newcommand{\bbz}{{\mathbb Z}}	
	
\newcommand{\bbq}{{\mathbb Q}}

\begin{document}
  \author{J\"{o}rg Hofmann}
  \address{Johannes Gutenberg-Universit\"{a}t Mainz, Institut f\"{u}r Mathematik\\ Staudingerweg 9, 55099 Mainz, Germany}
  \email{hofmannj@mathematik.uni-mainz.de}
  \title{Uniformizing differential equations of arithmetic $(1;e)$-groups}
  \begin{abstract}
   We present a numerical method to compute accessory parameters of the uniformizing differential equations of elliptic curves associated to arithmetic Fuchsian $(1;e)$ groups. Up to $PSL_2(\bbr)$ conjugation there are only finitely many such groups and we give a complete list of the accessory parameters.
  \end{abstract}
  \subjclass{34M35, 33E10 ,20H10}
  \keywords{Lam\'{e} equations, arithmetic Fuchsian groups, Uniformization}
  \date{\today}
  \maketitle

\section*{Introduction}
A discrete subgroup $\Gamma$ of $PSL_2(\bbr)$ is called \emph{Fuchsian group} and its elements act on the upper half plane $\bbh$ as M\"{o}bius transformations. The quotient $Y(\Gamma) = \bbh/\Gamma$ can always be equipped with a complex structure and is therefore a Riemann surface. An important class of Fuchsian groups are congruence subgroups of $SL_2(\bbz)$ which lead to modular curves, which can be compactified by adding finitely many cusps. An immediate generalization of congruence subgroups are \emph{arithmetic Fuchsian groups}. A group is called \emph{arithmetic} if it is commensurable with the embedding of the norm one elements of an order of a quaternion algebra into $M_2(\bbr)$. If $\Gamma$ is a an arithmetic Fuchsian group, then $X(\Gamma)$ can be realized as a projective algebraic curve defined by equations with coefficients in a number field. In section \ref{sec:shimura} we briefly review the basic facts on arithmetic Fuchsian groups. Special attention will be paid on the case where $\Gamma$ is a \emph{Fuchsian group of $(1;e)$-type}, that is a Fuchsian group whose quotient has genus one and the projection $\phi:\bbh \rightarrow \bbh/\Gamma$ branches exactly at one point and the branching index at this point equals $e$. Arithmetic Fuchsian groups of $(1;e)$ type were classified by Takeuchi in \cite{Tak83}. The main issue will be the inverse $\omega$ of $\phi$. As will be explained in section \ref{sec:unieq}, $\omega$ can be realized as the quotient of two linearly independent solutions of a differential equation $L$, called the \emph{uniformizing differential equation}. Recently Sijsling gave a list of equations for projective curves that yield models for most of the Riemann surfaces associated to the groups from Takeuchis list \cite{Sij2011}. If such a model is known, it determines the local data of the uniformizing equation. But in general, one complex parameter, called the \emph{accessory parameter}, remains undetermined. In section \ref{sec:num} we explain how a method that computes high precision approximations of the generators of the monodromy group of a linear Fuchsian differential equation can be used to approach the determination of the accessory parameter. In all cases we are able to compute floating point approximations of the accessory parameter with accuracy sufficiently high to identify the accessory parameters as algebraic numbers.

\section{Riemann Surfaces and Arithmetic Fuchsian groups of \texorpdfstring{$(1;e)$}{[(1;e)]}  type} \label{sec:shimura}
\begin{figure}[ht] 
  \includegraphics{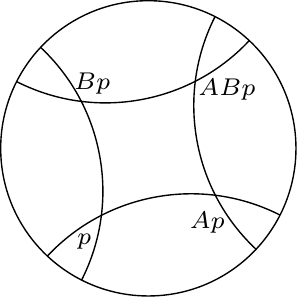}
\caption{Fundamental domain of a $(1;e)$ group} \label{fig:funda}
\end{figure}

One way to construct Riemann surfaces is by taking the quotient $\bbh/\Gamma$ by a suitable group acting on the upper half plane. We will focus on the case where $\Gamma$ is an arithmetic Fuchsian group of $(1;e)$-type. The books \cite{vign1980} and \cite{Katok} cover all the material needed.\\
We start by recalling the necessary definitions.
For a Fuchsian group $\Gamma$ whose action on $\bbh$ has a compact fundamental domain the quotient $\bbh / \Gamma$ is also compact with genus $g$. Since then necessarily $\Gamma$ has no parabolic elements, if the projection $\bbh/\Gamma \rightarrow \bbh$ branches above $r$ points with index $m_1,\ldots,m_r$ the \emph{signature} of $\Gamma$ can be defined as the tupel $(g;m_1,\ldots,m_r)$. Groups with signature $(0; m_1, m_2,m_3 )$ are called \emph{triangle groups} and their investigation has a long history. Beyond triangle groups a prominent role is played by groups with \emph{signature} $(1;e)$.
\begin{definition}[Fuchsian $(1;e)$-group]
A subgroup $\Gamma \leq PSL_2(\bbr)$ is called Fuchsian $(1;e)$-group if the quotient $\bbh/\Gamma$ of the action of $\Gamma$ on the complex upper half plane is a compact Riemann surface of genus one and the projection $\bbh \rightarrow \bbh/\Gamma$ ramifies of index $e$ above exactly one point.
\end{definition}
Following Fricke and Klein \cite{FrickeKlein}, a lift of a Fuchsian $(1;e)$-group $\Gamma$ to $SL_2(\bbr)$ can be presented as 
$$
   <\alpha,\beta,\gamma| \alpha \beta \alpha^{-1} \beta^{-1} \gamma =-id_2, \gamma^e =-id_2>, 
$$
 where $\alpha$ and $\beta$ are hyperbolic elements of $SL_2(\bbr)$ and $\gamma$ is elliptic with $tr(\gamma) = 2 \cos(\frac{\pi}{e})$. In addition the Fricke relation
$$
 tr(\alpha)^2 + tr(\beta)^2 + tr(\alpha \beta)^2 - tr(\alpha)tr(\beta)tr(\alpha \beta) = 2 - 2\cos\left(\frac{\pi}{e}\right)
$$
holds.
For realizations $A, B \in SL_2{\bbr}$ of $\alpha$ and $\beta$ and $p$ the fixed point of the elliptic element $ABA^{-1}B^{-1}$, the matrices $A$ and $B$ considered as maps are the side pairings of the hyperbolic polygon $P$ with vertices $p,Ap,Bp,ABp$ and angle sum $2\frac{\pi}{e}$. It is Poincar\'{e}'s theorem which in this case states that $P$ is a fundamental domain for $\Gamma$. After suitable choices the fundamental domain of a Fuchsian $(1;e)$ group in the upper half plane can be depicted as a quadrilateral, as in figure \ref{fig:funda}.
There are even up to conjugation infinitely many such groups, but if we impose the condition that the group is \emph{arithmetic}, a notion we will make precise below, only finitely many conjugation classes exist. The notion of \emph{quaternion algebras} is crucial. 
\begin{definition}[quaternion-algebra]
    For two elements $a,b$ of a field $F$ of characteristic not equal to $2$, the algebra $Fi \oplus Fj \oplus Fk \oplus Fl $ with multiplication given by $i^2=a,\ j^2=b$ and $ij=ji=k$ is denoted by $\genfrac(){}{}{a,b}{F}$ and called the  \emph{quaternion algebra} determined by $a$ and $b$. 
\end{definition}
Two familiar examples of quaternion algebras are the Hamiltonian quaternion algebra $\genfrac(){}{}{-1,-1}{\bbr} $ and the matrix algebra $\genfrac(){}{}{1,1}{F} \simeq M_2(F)$. The next lemma shows that after an extension to the real numbers, quaternion algebras are governed by the above two examples. 
\begin{lemma}
    Let $\phi_i$ be an embedding of the number field $F$ in $\bbr$, then $\genfrac(){}{}{a,b}{F} \otimes_{F,\phi_i} \bbr$ is either isomorphic to $M_2(\bbr)$ or a division algebra.
\end{lemma}
If $\genfrac(){}{}{a,b}{F} \otimes_{F,i} \bbr \cong M_2(\bbr)$ the embedding $\phi_i$ is named \emph{split} and \emph{nonsplit} elsewise.
In the following let $F$ be a totally real number field of degree $n$. There are $n$ distinct embeddings $\phi_i $ of $F$ into the real numbers, called infinite places.
Consider quaternion algebras which are split at the identity $\phi_1$ and non split at all other infinite places. Denote the corresponding isomorphism $\genfrac(){}{}{a,b}{F} \otimes_{F,1} \bbr \rightarrow M_2(\bbr)$ by $j$.
\begin{definition}[Order]
    The ring of integers of $F$ are denoted $Z_F$. We say that a $Z_F$ submodule of $\genfrac(){}{}{a,b}{F} \otimes_{F,i} \bbr$ is an \emph{order} if it is also a subring. 
\end{definition}Clearly an element $x$ of a quaternion algebra $\genfrac(){}{}{a,b}{F}$ can be written as $x =f_0i + f_1j + f_2k + f_3l $ and 
   $$|x| := f_0^2 - a f_1^2 -bf_2^2 + abf_3^2 $$
   introduces a norm. Denote all elements of norm $1$ of an order $\mathcal{O}$ by $\mathcal{O}^1$. 
   \begin{definition}[arithmetic Fuchsian group]
    A Fuchsian $(1;e)$-group $\Gamma$ is called arithmetic if there exists a quaternion algebra with order $\mathcal{O}$ such that $\Gamma$ is commensurable with $j(\mathcal{O}^1)$.   
   \end{definition}
   Takeuchi used that an arithmetic Fuchsian $(1;e)$-group $\Gamma$ can be described completely in terms of the triple 
   $$
     (tr(\alpha) , tr(\beta) , tr(\alpha \beta))
   $$ 
   to prove that up to $PSL_2(\bbr)$-conjugation there are 73 arithmetic Fuchsian groups of $(1;e)$-type.
   A set of generators $(\alpha_0,\beta_0,\gamma_0)$ of $\Gamma$ is called \emph{fundamental}, if 
   $$
    2< tr(\alpha_0) \leq tr(\beta_0) \leq tr(\alpha_0\beta_0) 
   $$
   and if it is of least height among all such generating triples, where the height is defined as 
   $$
    h(\alpha, \beta, \gamma) := tr(\alpha)^2 + tr(\beta^2) + tr(\alpha \beta)^2.
   $$ 
   In \cite{Tak83} a list of these fundamental triples $(x,y,z)$ for arithmetic Fuchsian $(1;e)$-groups is given.
   Moreover, the quotient $\bbh/\Gamma$ is an Riemann surface with a marked point, the point above which the projection $\bbh \rightarrow \bbh/\Gamma$ branches. Hence this quotient is an instance of an orbifold.

\section{Orbifold Uniformization and Differential Equations} \label{sec:unieq}
\begin{figure}[ht]
  \includegraphics{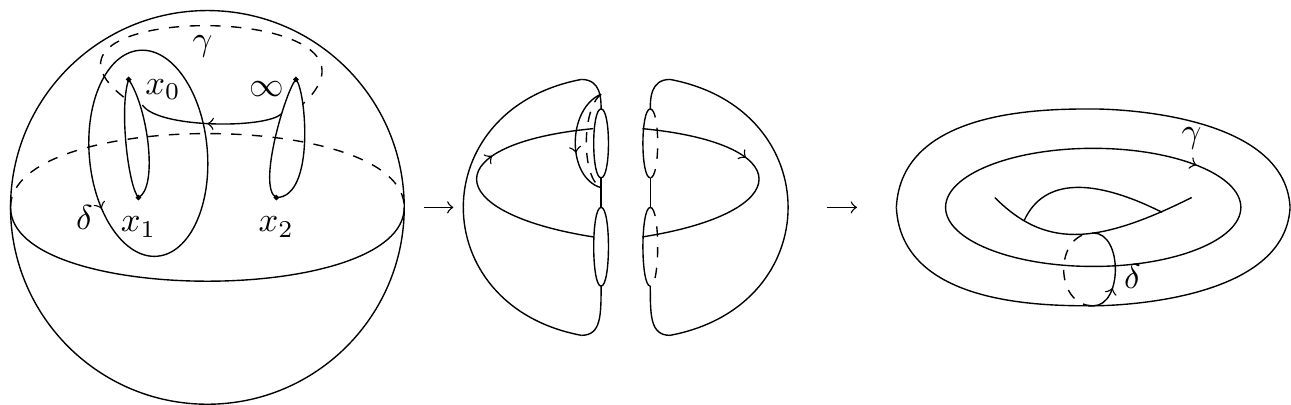}
\caption{Connection of $\pi_1(E)$ and $\pi_1(\bbp^1 \setminus \Sigma)$}
\label{fig:torusPi}
\end{figure}

In this section we recall shortly the classical theory of orbifold uniformization. A reference which explains the mathematical as well as the historical aspect of the uniformization of Riemann surfaces is the book \cite{Saint-Gervais} by the collective Henri Paul Saint Gervais.
\begin{definition}[orbifold]
Let $X$ be a complex manifold and $Y \subset X$ a hypersurface, that splits as $Y= \cup_j Y_j$ in irreducible components. Furthermore associate to every $Y_j$ a natural number $b_j \geq 2$ or $\infty$. The triple $(X,Y,(b_j)_j)$ is called an orbifold if for every point in $X \setminus \cup_j \left\{Y_j|b_j=\infty \right\}$ there is an open neighbourhood $U$ and a covering manifold which ramifies along $U \cap Y$ with branching indices given by $(b_j)_j$.
\end{definition}
If the above local coverings can be realized globally this definition can be subsumed as.
\begin{definition}[uniformization]
 If there exist a complex manifold $M$ and a map $\phi: M \rightarrow X$ which ramifies exactly along the hypersurfaces $Y_j$ with the given braching indices $b_j$, then the orbifold $(X,Y,(b_j)_j)$ is called uniformizable and $M$ is called a uniformization.
\end{definition}
The multivalued inverse of $\phi$ is called \emph{developing map} and it is in general a hard problem to determine this map explicitly. 
We will focus on the case where $X$ is the quotient of the upper half plane by an arithmetic Fuchsian group of $(1;e)$-type and therefore $M$ coincides with upper half plane itself. The hypersurface $Y$ is just one point above which the map $\phi$ ramifies with index $e$. 
There is a fundamental connection between the multivalued inverse of the map $\phi$ and a certain linear differential equation, which goes back to Poincar\"{e}. For a local coordinate $x$ on $X$ the inverse map $\omega=\phi^{-1}(x)$ is $PGL_2(\bbr)$-multivalued, but the Schwarzian derivative 
$$ 
(Sw):= \left( \frac{\omega^{'''}} {\omega^{'}} \right) - \frac{3}{2} \left( \frac{\omega^{''}} {\omega^{'}} \right)^2 
$$ 
is a single-valued map on $X$.
For a proof of the last two claims and the subsequent desired proposition connecting the Schwarzian derivative and differential equations see \cite{Yoshida1987}[Ch.4.2].
\begin{proposition}
Let $\omega(x)$ be a non-constant $PGL_2(\bbr)$-multivalued map, then there exist two $\bbc$-linearly independent solutions $y_0(x),y_1(x)$ of the differential equation 
\begin{align*}  
	L_E := \left[ \frac{d^2}{dx^2} + \frac{1}{2}S(\omega(x) )\right] y=0 , 
\end{align*}
such that $\omega(x) = \frac{y_1(x)}{y_0(x)}$.  
\end{proposition}
\begin{definition}[uniformizing differential equation]
 The differential equation from the above proposition is called uniformizing differential equation or Schwarzian differential equation.
\end{definition}
In general the unformizing differential equation can not be computed from $\Gamma$ directly, but one can try to recover $L_E$ from its monodromy group $M$, since the projectivization of $M$ coincides with $\Gamma$. 
To be more precise we recall the definition of the monodromy group of a linear homogeneous differential equation $L$.  Such an equation $L$ of order $n$ on a Riemann surface $S$ has near any ordinary or regular singular point $p$ exactly $n$ over $\bbc$ linear independent solutions, denote the space spanned by these solutions by $F_p$. Since analytic continuation of a germ of $\mathcal{O}_S$ along a path on $S$ is defined up to homotopy, there is a well defined representation 
$$
\rho:\pi_1(S \setminus \Sigma , p ) \rightarrow GL_n(\bbc),
$$
where $\Sigma$ is the set of singular points of $L$. 
For a loop $\gamma \in S$ with base point $p$, analytic continuation along $\gamma$ transforms $F_p$ into $\tilde{F}_p$ and since $F_p$ and $\tilde{F}_p$ are two bases of the same vector space, the solution space of $L$, analytic continuation yields an element  $M \in Gl_n(\bbc)$ as
$$
    F_p =  \tilde{F}_p M .
$$
\begin{definition}[monodromy representation, monodromy group]
The representation $\rho$ is called \emph{monodromy representation} of $L$ and  its image in $GL_n(\bbc)$ is called the \emph{monodromy group}. 
\end{definition}
Denote the set of branch points of the projection $\bbh \rightarrow X$ by $B$ and let $p$ a point outside of $B$, the representation $ \rho_0: \pi_1(X \setminus B ,p) \rightarrow \Gamma \subset PSL_2(\bbr)$ associated to the covering $\bbh \rightarrow X$ is conjugated to the projectivization of the monodromy representation of $\tilde{\rho}: \pi_1(X \setminus B,p) \rightarrow GL_2(\bbc) \rightarrow PGL_2(\bbc)$. 
This can be seen as follows. Two solutions $y_1, y_0$ of $L_E$ are a local branch of the inverse $\omega(x)$ of the covering map. By analytic continuation along a loop $\gamma$ centred at $p$ two $\bbc$-linear independent local  solutions $y_1,y_0$ of $L_E$ are converted to $\tilde{y_0} = ay_1 + by_0$, $\tilde{y_1} = cy_1 + dy_0$ with $M=\bigl( \begin{smallmatrix} a & b \\ c & d \end{smallmatrix} \bigr)\in GL_2(\bbc)$. This yields
$$
\frac{\tilde{y_1}}{\tilde{y_0}} = \frac{a y_1/y_0 + b}{c y_1/y_0 + d},
$$
which is again the quotient of two  solutions of $L_E$ and therefore another branch of $\omega(x)$ related to the former one by the M\"{o}bius transformation associated to $M$. Hence if the monodromy representation is computed according to an arbitrary basis $F_q$ the two representations $\rho_0$ and $\tilde{\rho}$ are conjugated.\\
We will do no computations with $L_E$ directly, instead we consider a differential equation $L_{\bbp^1}$ on the projective line which is related to $L_E$ via a certain pull-back.   
Since we assumed $X$ to be a compact Riemann surface it can be realized as an elliptic curve, that in an affine chart can be given in the form 
$$
      E:y^2=4x^3+ax+b,\ a,b \in \bbc.
$$
The map $\pi:E \rightarrow \bbp^1 \setminus \Sigma,\  (x,y) \mapsto x$  is a twofold cover of $\bbp^1 \setminus \left\{x_1,x_2,x_3,\infty \right\}$, branching of index two at the points of $\Sigma:=\left\{x_1,x_2,x_3,\infty \right\}$, where $x_i$  are the roots of $P(x):=4x^3+ax+b$. Thus $L_E$ is the pullback of a differential equation on $\bbp^1$, furnished with a local coordinate $t$ which is defined by local data as
\begin{align*}
   L_{\bbp^1}:=\left[ P(t)\frac{d^2}{dt^2} + \frac{1}{2} P(t)^{'} \frac{d}{dt} + n(n+1)t +C \right]y=0,\ n=\frac{1}{2e} - \frac{1}{2}.
\end{align*}
This differential equations is called \emph{algebraic Lam\'{e} equation} and the complex number $C$ is referred to as \emph{accessory parameter}. This kind of differential equation was first investigated in connection with ellipsoidal harmonics \cite{Poole}[ch. IX] .  The accessory parameter does not affect the local exponents and the location of the singularities of $L_{\bbp^1}$. That $L_{\bbp^1}$ is indeed the correct object can be explained by a closer look at its Riemann scheme 
    \begin{align*}
      \begin{Bmatrix} t_1 & t_2 & t_3 & \infty \\ \hline
                      0   & 0   & 0   & -\frac{n}{2} \\
                      \frac{1}{2}   & \frac{1}{2}   & \frac{1}{2}   & \frac{n+1}{2} \\
      \end{Bmatrix}.
    \end{align*}
A basis of the solution space of $L_{\bbp^1}$ at a regular point can be given by two power series and in a vicinity of any of the finite singular points it is given by 
    \begin{align*}y_0(t) &= \sum_{n=0}^{\infty} a_n (t-t_i)^n \ \textrm{and} \\
                  y_1(t) &= (t-t_i)^{\frac{1}{2}}\sum_{n=0}^{\infty} b_n (t-t_i)^n,
    \end{align*}
whose quotient is $\frac{y_1(t)}{y_0(t)} = (t-t_i)^{\frac{1}{2}} \sum_{n=0}^{\infty} c_n (t-t_0)^n$ \cite{ince}. Similarly since $-\frac{n}{2} - \frac{n-1}{2} = -\frac{1}{2e}$ the shape of this quotient at infinity is 
$$
    \left(\frac{1}{t}\right)^{\frac{1}{2e}}\sum_{n=0}^{\infty} c_n \left( \frac{1}{t}\right)^n.
$$
At an ordinary point $p$ of $L_{\bbp^1}$ the local exponents are $0$ and $1$, hence following \cite{ince}[p.396 ff] in a vicinity of $p$ the Ansatz 
   $$
    y_{\sigma}(t)  = \sum_{i=0} ^{\infty}  a_n(\sigma) (t-p)^{n+\sigma}=0,\ \sigma \in \left\{ 0,1\right\} 
   $$ 
provides a basis of the solution space of $L_{\bbp^1}$ consisting of convergent power series, adding derivatives a fundamental matrix can be found as 
$$
F_p = \begin{pmatrix} 
			y_1(t) & y^{'}_1(t) \\ y_0(t) & y^{'}_0(t) 
      \end{pmatrix}.
$$
The radius of convergence equals $\bigl|\Sigma,p \bigr|$, the minimal distance from $p$ to any of the singular points of $L_{\bbp^1}$. The coefficients are computed recursively as 
\begin{align*}
   a_n(\sigma) =-&\frac{ g_2( (j-1+\sigma)(j-2+\sigma)+(j-1+\sigma) )a_{j-1}}{g_3(n+\sigma)(n-1+\sigma)} \\+ 
         &\frac{( n(n+1) + C - 4(j-2+\sigma)(j-3+\sigma))a_{j-2}}{g_3(n+\sigma)(n-1+\sigma)} \\-
         &\frac{(6(j-3+\sigma))a_{j-3} }{g_3(n+\sigma)(n-1+\sigma)}, \ n \geq 3.
\end{align*}
where the initial conditions, i.e. the values of $a_0, a_1$ and $a_2$ could be chosen arbitrarily with at least one $a_i\neq0$, but we will fix them as $a_0=1, a_1=0$ and $a_2=0$.  
The monodromy group of $L_{\bbp^1}$ is generated by the images of the monodromy presentation of loops $\gamma_1,\gamma_2,\gamma_3$ starting at $p$ and encircling exactly one of the finite points in $\Sigma$ in counterclockwise direction, with the additional property that composition of pathes $\gamma_1 \gamma_2 \gamma_3 $ is homotopic to a path $\gamma_{\infty}$ encircling $\infty$ once in clockwise direction. The images of theses pathes under $\rho$ are called $M_1,M_2,M_3$ and $M_{\infty}$.
To compute the analytic continuation choose points $p_0,\ldots,p_{m-1},p_m=p_0$ on (a path homotopic to) $\gamma$ s.t. $p_{i+1}$ lies in the radius of convergence of the elements of $F_{p_i}$. With the above choice of initial conditions $F_{p_i}(p_i)=id_2$ holds and
   $$
    F_{p_i}(t) = F_{p_{i+1}}(t) F_{p_i}(p_{i+1}), 
   $$  
is implied. Which finally leads to  
$$
  M_{\gamma} = \prod_{i=0}^{m-1}  F_{p_{m-i}}(p_{m-i-1}).
$$ 
To relate the monodromy representation of $L_{\bbp^1}$ with our initial group $\Gamma$ it is enough to understand how the monodromy group of $\bbp^1  \setminus \Sigma$ lifts to the fundamental group of $X$.
We recall the standard construction of the fundamental group of an pointed elliptic curve as branched two fold covering of the four times punctured sphere. Cut the Riemann sphere from $t_1$ to $t_2$ and from $t_2$ to $\infty$, expand the cuts a little and glue two copies of this cutted sphere along the cuts with opposite orientation. The fundamental group of the pointed elliptic curve is generated by lifts of the loops $\delta$ and $\gamma$ as in figure \ref{fig:torusPi}. Hence the monodromy of $L_E$ along $\delta$ coincides with the monodromy of $L_{\bbp^1}$ along a loop $\gamma_1 \gamma_2$ and the monodromy of $L_E$ along $\gamma$ coincides with the monodromy of $L_{\bbp^1}$ along $ \gamma_{\infty} \gamma_{2}^{-1} \sim \gamma_3 \circ \gamma_1$. Thus we can identify up to scalars $\Gamma$ with the group generated by $M_3 M_1$ and $M_1 M_2$. Note that $M_3 M_1 M_1 M_2 = M_3 M_2$. According to Fricke the traces of $M_i,i=1,\ldots,3$ and the traces of the products above determine the monodromy representation of $L_{\bbp^1}$.
This relation between the monodromy of $L_{\bbp^1}$ and the Fuchsian group $\Gamma$ can also been established, in the spirit of Klein and Poincar\'{e}, by a look at the image of a quotient $y = \frac{y_1}{y_0}$ of the upper half. Then suitable analytic continuation of $y$ and the Schwarz reflection principle would exhibit how certain elements of the monodromy group of $L_{\bbp^1}$ coincide with side pairings of the fundamental domain of $\Gamma$. This is exactly the way how the monodromy of the hypergeometric function is related to triangle groups. 
Four points on the Riemann sphere can always be mapped to $0,1,A,\infty$ by a M\"{o}bius transformation, that means given a Fuchsian $(1;e)$-group $L_{\bbp^1}$ depends on two parameters, namely $A$ and the accessory parameter. It is explained in \cite{Ihara} that $A$ and $C$ have to be elements of $\overline{\bbq}$ and moreover in a suitable coordinate they have to be elements of a number field that can be determined from $\Gamma$ directly. But recently J. Sijsling \cite{Sij2011} constructed Weierstra\ss -equations for members of isomorphism or at least isogeny classes of the elliptic curves associated to those 73 groups which provides the local data of the uniformizing differential equation i.e. the location of the singular points and the Riemann scheme of the corresponding Lam\'{e} equations are known in these cases.
Whenever $\Gamma$ is commensurable with a triangle group he additionally used in \cite{Sij2011i}  Bely\u{i} maps to determine the accessory parameters of the uniformizing differential equations. We will use his results and numerical methods in the next section to tackle the problem of the determination of accessory parameters by investigating the monodromy of the uniformizing differential equations.

\section{Approximation of Accessory Parameters} \label{sec:num}
In \cite{chud88} G.V. and D.V. Chudnovsky considered algorithms to compute approximations of generators of the monodromy group of linear differential equations. They also had available implementations and used them to do high precision approximations to guess accessory parameters in at least a few cases picked from the Lam\'{e} equations associated to Takeuchis list \cite{chud89}. In a talk in Banff in 2010 \cite{Beu2010} F. Beukers stated that Sijsling used similar numerical methods and was again successful in a few cases. But to the authors knowledge there is no complete list of accessory parameters available. Hence we picked up that approach and did our own implementations as follows. Once we gained the insights from the last section the approach is straight forward. Replace
the solutions $y_1$ and $y_0$ of 
$$
L_{\bbp^1}:=\left[ P(t)\frac{d^2}{dt^2} + \frac{1}{2} P(t)^{'} \frac{d}{dt} + n(n+1)t +C \right]y=0,\ n=\frac{1}{2e} - \frac{1}{2}
$$
by truncations
$$
 y_{\sigma}^N(t)  = \sum_{i=0}^{N}  a_n(\sigma) (t-p)^{n+\sigma}=0,\ \sigma \in \left\{ 0,1\right\}
$$
and the fundamental matrix $F_p$ by
$$
F_p^N = \begin{pmatrix} 
			y_1^N(t) & {y_1^{N}}^{'}(t) \\ y_0^N(t) & {y_0^{N}}^{'}(t) 
      \end{pmatrix}.
$$
and approximate the monodromy along the path $\gamma$ as
$$
 A_{\gamma} = \prod_{i=0}^{m-1}  F_{p_{m-i}}^N(p_{m-i-1}).
$$
How this approximation is done effectively was investigated by G.V and D.V. Chudnovskys in \cite{chud87-2} and \cite{chud88}. There they considered and partially answered questions dealing with fast computations of the coefficients $a_n(\sigma)$ and how the points $p_i$ should be chosen. For details on the choices of loops representing a homotopy class and vertices of the involved polygons and only the implementations done see \cite{Hof1}. The errors occurring due to the truncation process can be bounded by explicit bounds.
As explained above the monodromy group of $L_{E}$ is generated by the products $M_1 M_2$ and $M_3 M_2$. In most cases the local data of $L_{\bbp^1}$ is known by Sijslings work, and the determination of the accessory parameter is left. It has to be determined such that the trace triple 
    \begin{align} \label{tracetr}
     	( tr(M_3 M_2),tr(M_1 M_2),tr(M_3 M_1) )  
    \end{align}
or a permutation hereof coincides with one of the triples in Takeuchis list. If only a permutation of the trace triple is found, one can use a M\"{o}bius transformation interchanging the singular points of $L_{\bbp^1}$ to recover the actual triple. If only the isogeny class of $E$ is known, we started with one member of this class and computed all elliptic curves $\tilde{E}$, such that there is an isogeny $\rho:E \mapsto \tilde{E}$ of degree bounded by a given $n \in \bbn$ and tried to find algebraic accessory parameters for differential equations with the local data specified by $\tilde{E}$. The tools provided by \cite{Sij-tools} proved to be very helpful.
Denote the tuple of generators of the monodromy of $L_{\bbp^1}$ with accessory parameter $C$ by $M(C)= (M_1(C),M_2(C),M_3(C))$ and its approximations by $A(C)= (A_1(C),A_2(C),A_3(C))$.  Since the monodromy group for the correct value of $C$ is a subgroup of $SL_2(\bbr)$, it is reasonable to consider the maps
\begin{align*}
     f(T)&:= max\left\{ Im(tr(T_1 T_2)) , Im(tr(T_3 T_1)) , Im(tr(T_3 T_2)) \right\} \quad   \text{and} \\
     g(T)&:= \|Im(tr(T_1 T_2))\| + \|Im(tr(T_3 T_1))\| + \|Im(tr(T_3 T_2))\| \\ 
\end{align*}
for a 3-tuple $T$ of $GL_2(\bbc)$ matrices.
By a theorem of Hilb \cite{Hilb}, there is an infinite but discrete subset of accessory parameters in the complex plane such that $f(M(C))=0$ and $g(M(C))=0$, 
but luckily numerical investigations suggest, that only the one which has smallest absolute value is algebraic.
If the coefficients of the elliptic curve $E$ are real, we use the function $f$ to approximate $C$. All computations are done in multiprecision arithmetic as provided by the $C$ library \emph{mpc} \cite{mpc} or the computer algebra system \verb!maple!. Start by choosing a small number $ex$ and $C_1$, such that the signs of $f(N(0))$ and $f(N(C_1))$ differ and build the pair $S=(S1,S2) = (f(N(0)) ,f(N(C_1)))$ which is used as input of algorithm \ref{alg1}.\\
\SetKw{Return}{return}
\SetKwInOut{Input}{input}
\SetKwInOut{ElIf}{elseif}
\begin{algorithm} \nllabel{alg1}
\caption{Algorithm for real coefficients}
   \Input{$S$}
   \BlankLine
   \uIf{$|f(N(C_1)|<ex$}{ 
    \Return $C_1$
   }
    \uElseIf{$|f(N(0)|<ex$ }{
     \Return 0
   } 
   \Else{ 
        $t= C_1/2$\\
        \While{$|f(N(t)|>ex$}{ 
          $t=(S1+S2)/2$\\
          \uIf{$sgn(f(N(t))) = sgn(f(N(S1)))$}{ 
            S=(t,S2)
 	 }
          \Else{
            S=(S1,t)
         }
         adjust N and the length of the mantissa
     }
      \Return t\\
    }
\end{algorithm}
If the accessory parameter is expected not to be real, use the map $g$. Choose a small real number $r$ and $C=0$ to compute $g(N(C))$ and continue with the code from algorithm \ref{alg2}.\\
\begin{algorithm} \nllabel{alg2}
\Input{t,C}
 \uIf{$g(N(C))<0$}{
  \Return{ $C$};
 }
 \Else{
  \While{$g(N(C))>ex$} {
    $cont=false;$ \\
   \For{$k=0$ to $7$}{
     $CC = C + r \exp(2\pi i/k)$\\
    \If{$g(N(CC) < g(N(C))$}{
     $C = CC;$\\ 
     $cont =true;$\\
     adjust $N$ and the length of the mantissa \\
     $ break;$ \\
     }
     \If{$cont=false$}{
      $ r=r/2; $\\ 
     }
   }
  }
   \Return{ $C$};   
 }
 \caption{Algorithm for complex coefficients}
\end{algorithm}
Adjust $N$ and the mantissa means that in every step the number of coefficients in the power series expansion of the solutions and the length of the mantissa are chosen as short as possible to guarantee that the errors made by the truncation and by cutting off the involved floating point numbers is smaller than the precision we have already achieved in the traces of $A(C)$ compared to the target traces $M(C)$.
The precision one needs to identify a floating point number as an algebraic number depends on the \emph{height} and the degree of its \emph{minimal} polynomial. An algorithms that accomplishes the identification of floating point numbers as algebraic expressions is described in \cite{PSLQ} and there exist various implementations for example in \verb!maple!.
Therefore the running time of the algorithms \ref{alg1} and \ref{alg2} for given $ex$ depends mainly on two factors. Namely on the height of the accessory parameter and in addition on the location of the singular points of $L_{\bbp^1}$. The influence of the first factor is obvious. The location of the singular points affects the running time as the configuration of the singular points mainly determines the number of expansion points $p_i$ needed in the approximation of the monodromy as explained above. Hence in the simplest cases after only a view minutes we obtained a promising candidate for the accessory parameter, but in the more complex cases the computer had to run for several hours. Once we obtained such a candidate $C_{alg}$ we checked that $A(C_{alg})$ and the trace triple from Takeuchis list coincid up to permutation up to at least 200 Digits. The reason why we can not turn our method in a rigorouse proof is that we do not have any a priori knowledge on the heights of the coefficients of $L_{\bbp^1}$.

The elliptic curves and accessory parameters found in this way are listed in the tables below. The labels are chosen according to \cite{Sij2011}. The symbol 
$$
n_d / n_D r
$$ 
encodes
\begin{itemize}
 \item $n_d$: the discriminant of the coefficient field $F$ of $\genfrac(){}{}{a,b}{F}$  
 \item $n_D$: the norm of the discriminant of the associated $\genfrac(){}{}{a,b}{F}$ 
 \item $r$: roman number used to distinguish cases with equal $n_d$ and $n_D$.
\end{itemize}
The data in the tables below is the label $n_d / n_D r$ an equation for the corresponding ellitic curve $E:\ y^2 + x^3 + a_1xy +a_3 y = x^3 +a_2x^2 +a_4x +a_6$, the accessory parameter $C$ of the pulled back uniformizing differential equation $L_{\bbp^1}$ and the the squares of the entries of the trace triple \eqref{tracetr}. The polynomial $P$ which specifies the coefficients of the differential equation $L_{\bbp^1}$ can be recovered from $E$ by the substitution of $y$ by $\frac{1}{2}(y-a_1x-a_3)$ and the elimination of the $x^2$ term. If an algebraic number $\gamma$ is involved in the coefficients of $E$ or the given trace triples has a complicated radical expression, we give some digits of its floating point expansion and its minimal polynomial $f_{\gamma}$. 
Whenever all nonzero singular points of $L_{\bbp^1}$ have non-vanishing real part the basepoint $p$ is chosen as the imaginary unit $i$, elsewise it is chosen as $0$. The singular points $t_i,\ i=1,\ldots,3$, which are the roots of $P$ are ordered according to the argument of the complex number $t_i-p$. This fixes the loops $\gamma_i$ and hence the trace triple $\eqref{tracetr}$.

\renewcommand{\arraystretch}{1.01}
\setlength{\extrarowheight}{2mm}
\begin{longtable}[l]{l|l l } 
\caption{ramification index 2}\\ 
\hline
\multirow{2}{*} {1/6i} & $ {y}^{2}+xy+y={x}^{3}-334x-2368 $ \\ 
  & $C=-{\frac {79}{64}}$   \\
  & $( 5,12,15 )$ \\[0.4mm]
\hline  
 
\multirow{2}{*} { 1/6ii }  &  $ {y}^{2}={x}^{3}-{x}^{2}-4x+4$ \\  
  & $ C=\frac{1}{16} $    \\
  &  $\left(8,6,12 \right)$  \\[0.6mm] 
\hline  
 
\multirow{2}{*} { 1/14 }  &   $  y^2=x^3+23220x-2285712$ \\  
  & $  C= {\frac {9}{8}}$    \\
  & $ \left(7,7,9\right) $\\[0.6mm]  
\hline  

\multirow{3}{*} {5/4i}  &  $ {y}^{2}={x}^{3}- \left( -\frac{1}{2}\sqrt{5}+\frac{1}{2}\right) {x}^{2}- \left(8+3\sqrt{5}\right) x-{\frac{25}{2}}-\frac{11}{2}\sqrt {5}$ \\  
  & $  C=-{\frac {5}{32}}-\frac{1}{32}\sqrt{5} $    \\ 
  &$  \left(3+\sqrt{5} ,6+2\sqrt{5}, 7+3\sqrt{5} \right) $ \\[0.6mm]  
\hline  

\multirow{2}{*} {5/4ii}  &  ${y}^{2} + \left( \frac{3}{2}+\frac{1}{2}\sqrt {5}\right) y={x}^{3}- \left( -\frac{1}{2}\sqrt{5}+ \frac{1}{2}\right) {x}^{2}-\left(-\frac{329}{2}+\frac{111}{2}\sqrt {5} \right) x +\frac{769}{2}-\frac{287}{2}\sqrt{5}$ \\  
  & $  C= -\frac{5}{128}+\frac{1}{128}\sqrt {5}$  \\
  & $ \left( \frac{7}{2}+\frac{3}{2}\sqrt{5} ,6+2\sqrt{5}, \frac{7}{2}+\frac{3}{2}\sqrt{5}\right) $  \\[0.6mm]  
\hline

\multirow{2}{*} { 5/4iii }  & $  {y}^{2}+xy+y={x}^{3}-\left(-\frac{1}{2}\sqrt{5}+\frac{1}{2}\right) {x}^{2}- \left(\frac {101}{2}+\frac{45}{2}\sqrt{5}\right)x-{\frac{1895}{2}}-{\frac{847}{2}}\sqrt{5}$ \\  
& $  C=-\frac{5}{64}-\frac{1}{32}\sqrt{5}$  \\ 
& $ \left(\frac{9}{2}+\frac{3}{2}\sqrt {5} ,\frac{9}{2}+\frac{3}{2}\sqrt {5},\frac{7}{2}+\frac{3}{2}\sqrt{5} \right) $ \\[0.6mm]  
\hline

\multirow{2}{*} { 8/7i }  &  $ {y}^{2}={x}^{3}+\sqrt{2}{x}^{2}- \left( 142\sqrt{2}+202 \right) x-1170\sqrt{2}-1655 $  \\  
& $ C=-\frac{15}{16}-\frac{5}{8} \sqrt{2} $   \\ 
& $\left(3+\sqrt{2} ,12+8\sqrt {2},13+9\sqrt {2} \right)$ \\[0.6mm]  
\hline  
  
\multirow{2}{*} { 8/7ii }  & $ {y}^{2}={x}^{3}+\sqrt{2}{x}^{2}- \left( 142\sqrt{2}+202 \right) x-1170\sqrt{2}-1655 $  \\  
& $ C=-\frac{15}{16}+\frac{5}{8} \sqrt{2} $   \\ 
& $\left( 3+2\sqrt {2} ,5+3\sqrt {2}, 5+3\sqrt {2} \right)$ \\[0.6mm] 
\hline 

\multirow{3}{*} {  8/2 }  &  \multicolumn{2}{l}{  $  {y}^{2}={x}^{3}-x$ }\\  
& $C = 0 $  &  \\ 
& $ \left(3+2\sqrt {2}  ,5+3\sqrt{2},5+3\sqrt{2} \right)$&\\[0.6mm]  
\hline

\multirow{2}{*} {12/3}  & \multicolumn{2}{l}{$ {y}^{2}={x}^{3}- \left( \sqrt {3}-1 \right) {x}^{2}- \left( 65\sqrt {3}+111 \right) x+348\sqrt {3}+603$ }\\  
& $ C=\frac{5}{8}+\frac{5}{16}\sqrt {3}$  &   \\
& $\left(8+4\sqrt{3} ,3+\sqrt {3},9+5\sqrt{3} \right) $ \\[0.6mm]  
\hline  

\multirow{2}{*} { 12/2 } &  \multicolumn{2}{l}{  $  {y}^{2}={x}^{3}+a{x}^{2}+x+3\,\sqrt {3}-5$  }\\  
& $C=0$  & \\
& $ \left(4+2\sqrt {3} ,4+2\sqrt {3}, 4+2\sqrt {3}\right) $   \\[0.6mm]  
\hline

\multirow{3}{*} {13/36} & \multicolumn{2}{l}{  $ {y}^{2}+ \left( \frac{3 +\sqrt {13}}{2} \right) y={x}^{3}-\frac{5601845}{2}- \left(
\frac{60077-16383\sqrt{13}}{2} \right) x+\frac{1551027}{2}\sqrt{13}$ }\\  
& $ C= -{\frac {1625}{128}}+{\frac {375}{128}}\,\sqrt {13}$ &  \\[0.1cm]
& $\left(\frac{5}{2}+\frac{1}{2}\sqrt {13} ,16+4\sqrt{13},\frac{33}{2}+\frac{9}{2}\sqrt{13}\right)$ & \\[0.6mm]

\hline  

\multirow{2}{*}{13/4}  & $ {y}^{2}+\left(\frac{3+\sqrt{13}}{2}\right)y={x}^{3}-\left(\frac{3+\sqrt {13}}{2}\right){x}^{2}-\left(\frac{275+75\sqrt{13}}{2} \right) x-\frac{1565-433\sqrt {13}}{2}$ \\  
& $C=-\frac{65}{128}-\frac{15}{128}\sqrt{13}$    \\  
&  $\left(\frac{11+3\sqrt {13}}{2} ,\frac{11+3\sqrt {13}}{2},\frac{7+\sqrt {13}}{2} \right)$\\[0.6mm]
\hline

\multirow{2}{*} {17/2i}  &  $ {y}^{2}+xy+ \left( \gamma+1 \right) y={x}^{3}+a{x}^{2}- \left( -61\gamma+157 \right) x+348\gamma-896$ \\
& $\gamma = \frac{1-\sqrt{17}}{2}$  \\
& $C=\frac{-55+20\gamma}{64} $  \\
& $ \left(\frac{5+\sqrt{17}}{2},10+2\sqrt{17},\frac{21+5\sqrt{17}}{2}\right)$\\[0.6mm]  
\hline

\multirow{2}{*} {17/2ii} & $ {y}^{2}+xy+ \left( a+1 \right) y={x}^{3}+a{x}^{2}- \left( -61\gamma+157 \right) x+348\gamma-896 $ \\  
& $\gamma = \frac{1+\sqrt{17}}{2}$  \\
& $C=\frac{-55+20\gamma}{64} $ \\  
&  $\left(\frac{7+\sqrt{17}}{2},5+\sqrt{17},\frac{13+3\sqrt{17}}{2} \right)$\\[0.6mm]
\hline  

\multirow{2}{*} {21/4}  & 
$ {y}^{2}+ \left( \frac{ \gamma^2+19\gamma+15}{14} \right) y=x^3-\left( \frac{ \gamma^2-9\gamma+15}{14} \right) x^2 
- \left( \frac{99\gamma^3-144\gamma^2-618\gamma-1451} {14}\right) x $ \\
&$-\frac{88\gamma^3-505\gamma^2-777\gamma-555}{14}$\\
& $f_{\gamma}  =  x^4-4x^3-x^2+10x+43, \quad \gamma = - 1.292- 1.323i$ \\
& $  C=\frac{15 - 5\sqrt {21} + 56 i\sqrt{3} - 24 i \sqrt {7}}{128} $ \\  
&$ \left( \frac{15+3\sqrt{21}}{2} ,\frac{5+\sqrt{7}\sqrt{3}}{2},\frac{17+3\sqrt{7}\sqrt{3}}{2} \right)$\\[0.6mm]
\hline

\multirow{2}{*}{24/3} & $ y^2=x^3-\left( -\frac{1}{16}\gamma^{3}-\frac{1}{8}\gamma^{2}-\frac{1}{2}\gamma+1 \right) x^2 - \left( \frac {17295}{4}\gamma^3-\frac{14243}{8}\gamma^2-60459\gamma+155218 \right) x$ \\
& $+\frac{8148639}{8}\gamma^{3}-697026\gamma^2-\frac{25765125}{2}\gamma+34343808 $\\
& $f_{\gamma} = x^4-8x^2+64$,  $\gamma = - 2.450- 1.414i$ \\
& $C= \frac{4027}{288}-\frac{137}{24}\sqrt {6}+\frac{17}{36}i\sqrt{2}-\frac{37}{96}i\sqrt {3}$  \\
& $\left(5+2\sqrt{2}\sqrt {3},3+\sqrt{6},9+3\sqrt{6} \right)  $ \\[0.6mm]  
\hline

\multirow{2}{*} {33/12} & $ {y}^{2}+xy={x}^{3}- \left( \frac{-3-\sqrt{33}}{2} \right) x^2- \left( \frac{141}{2}
-\frac{27}{2}\sqrt{33} \right) x+\frac{369}{2}-\frac{63}{2}\sqrt{33}$  \\  
& $C=\frac{55}{576}-\frac{5}{288}\sqrt{33} $ \\ 
& $\left(\frac{7+\sqrt{11}\sqrt{3}}{2},\frac{9+\sqrt{33}}{2},  6+\sqrt{33}\right)$ \\[0.6mm]  
\hline

\multirow{2}{*} {49/56} &  $ y^2+xy+y=x^3-874-171x$ \\  
& $C=-{\frac {55}{64}}$ \\ 
& $\rho = 2 \cos(\frac{\pi}{7})  $ \\
& $\left(3\rho^2 + 2\rho-1  , 3\rho^2 + 2\rho-1 ,\rho^2 + \rho  \right) $ \\[0.6mm] 
\hline

\multirow{2}{*} {81/1} & $ y^2+xy+y=x^3-x^2-95x-697$ \\  
& $C = -\frac {15}{64}$ \\
& $f_{\rho} = x^3-3x-1, \quad \rho = 1.879$ \\
& $\left( \rho^2 + \rho +1    , (\rho+1)^2 , (\rho+1)^2\right)$ \\[0.6mm]  
\hline

\multirow{2}{*} { 148/1i} & $ y^2=x^3- \left( 464\gamma^{2}-320\gamma-1490 \right) x^2+x$ \\  
& $f_{\gamma} = x^3-x^2-3x+1,\quad \gamma = 2.170 $ \\
& $ C=\frac{1363+292\gamma-424\gamma^2}{16}$ \\
&$ \left( \gamma^2+\gamma, \gamma^2+2\gamma+1  ,\gamma^2+\gamma   \right)$ \\[0.6mm]
\hline  

\multirow{2}{*} {148/1ii} & $ y^2=x^3- \left( 464\gamma^{2}-320\gamma-1490 \right) x^2+x$ \\  
& $f_{\gamma} = x^3-x^2-3x+1, \quad  \gamma =  0.311$ \\
& $ C=\frac{1363+292\gamma-424\gamma^2}{16}$ \\  
& $ \left(-12\gamma^2+8\gamma+40,  -\gamma^2+\gamma+4,-13\gamma^2+9\gamma+42 \right)$ \\[0.6mm]
\hline

  \multirow{2}{*} { 148/1iii} & $ y^2=x^3- \left( 464\gamma^{2}-320\gamma-1490 \right) x^2+x$ \\  
  & $f_{\gamma} = x^3-x^2-3x+1, \quad \gamma = 1.481$ \\
  & $ C=\frac{1363+292\gamma-424\gamma^2}{16}$ \\  
  & $ \left( \gamma^2-3\gamma+2 , \gamma^2- 2 \gamma + 1 , \gamma^2-3\gamma+2 \right)$ \\[0.6mm]
  \hline  

\multirow{2}{*} {229/8i}  & $y^2=x^3-\left( 663\gamma^{2}-219\gamma-2485 \right) x^2- \left( -30778\gamma^2+13227\gamma+109691 \right) x$ \\  
  & $ C = \frac{205}{2}+\frac{105}{16}\gamma-\frac{105}{4}\gamma^2 $ \\
  & $f_{\gamma} = x^3-4x-1, \quad \gamma = 2.115$ \\ 
  & $\left( \gamma +2 , 8\gamma^2+16\gamma+4, 8\gamma^2+17\gamma+4 \right)$  \\[0.6mm] 
  \hline     

\multirow{2}{*} {229/8ii}  & $y^2=x^3-\left( 663\gamma^2-219\gamma-2485 \right) x^2- \left( -30778\gamma^{2}+13227\gamma+109691 \right) x$ \\  
  & $f_{\gamma} = x^3-4x-1, \quad \gamma = -0.254$ \\ 
  & $ C = \frac{205}{2}+\frac{105}{16}\gamma-\frac{105}{4}\gamma^2 $ \\  
  & $ \left(-3\gamma^2 + \gamma +13 , -\gamma^2 +5 , -4\gamma^2 + \gamma +16  \right)$  \\[0.6mm]  
  
\hline     

\multirow{2}{*} {229/8iii}  & $y^2=x^3-\left( 663\gamma^2-219\gamma-2485 \right) {x}^{2}- \left( -30778\gamma^2+13227\gamma+109691 \right) x$ \\  
  & $f_{\gamma} = x^3 - 4x -1, \quad \gamma = -1.860 $ \\
  & $ C =\frac{205}{2}+{\frac {105}{16}}\gamma-{\frac {105}{4}}\gamma^2 $ \\ 
  & $\left( \gamma^2 - 2 \gamma , \gamma^2 - 2\gamma +1, \gamma^2 - 2 \gamma    \right)$  \\[0.6mm]  
  
\hline

\multirow{2}{*} {725/16i} &   $  y^2+xy+\gamma y=x^3+x^2-\left( -447\gamma+4152 \right) x-85116\gamma + 59004$  \\
  & $f_{\gamma} = x^2-x-1, \quad \gamma =1.618 $ \\
  & $ C =\frac{ 205-300\gamma}{64}  $ \\
  & $f_{\rho} = x^4 - x^3 -3x^2 +x +1,\ \quad \rho = -1.355 $ \\ 
  & $\left( -2\rho^3 + 5\rho^2 -\rho - 1 , -2\rho^3 + 5\rho^2 -\rho - 1 , -\rho^3 +3\rho^2+\rho\right)$ \\[0.6mm]

 \hline

\multirow{2}{*} {725/16ii} & $ y^2+xy+\gamma y=x^3+x^2-\left( -447\gamma+4152 \right) x-85116\gamma + 59004$  \\
  & $f_{\gamma} = x^2-x-1, \quad \gamma =-0.618 $ \\
  & $ C =\frac{ 205-300\gamma}{64}  $ \\
  & $f_{\rho} = x^4 - x^3 -3x^2 +x +1,\ \quad \rho = -0.477 $ \\ 
  &$\left(9\rho^3-13\rho^2-21\rho+19 ,\rho^3 - 2\rho^2 - 2\rho + 4,9\rho^3-13\rho^2-21\rho+19 \right)$ \\[0.6mm]  
  
 \hline  

\multirow{2}{*} {1125/16} & ${y}^{2}+1/15  \left( \gamma^5+\gamma^4+17\gamma^3+14\gamma^2+13\gamma+19 \right) xy$ \\
&$+ \frac{1}{4629075} \left(\gamma^7+211776\gamma^6+471599\gamma^5+182985\gamma^4+3251185\gamma^3+8290968\gamma^2\right.$ \\
&$\left. +9151653\gamma+7962897 \right)y -  x^3-\frac{1}{4629075} \left(\gamma^7+520381\gamma^6+162994\gamma^5-125620\gamma^4\right.$\\
&$ \left. -1995100\gamma^3+3970498\gamma^2-3192547\gamma+5185452 \right) {x}^{2}$ \\
&$-\frac{1}{149325}\left(-12692863\gamma^7+86428787\gamma^6-164116067\gamma^5+518100715\gamma^4\right.$ \\                                                                      &$\left.-967426690\gamma^{3}+3757504646\gamma^{2}-3892822254\gamma+10486471269 \right)x$ \\
&$+\frac{10639239397}{925815}\gamma^7-\frac{26847806027}{925815}\gamma^6+\frac{17109768152}{185163}\gamma^5-\frac{4243037037}{28055}\gamma^{4
     }+\frac{595858143338}{925815}\gamma^3$ \\ 
   & $- \frac{316402671731}{308605}\gamma^2+\frac{2018798518646}{925815}\gamma-\frac{2030456532221}{925815} $ \\  

  &$ C=-\frac {2129}{1234420}a^7-\frac{42497}{29626080}a^6-\frac{50863}{7406520}a^5-\frac{173}{16833}a^4$\\&$-\frac{746393}{11850432}a^3-\frac{9461581}{59252160}a^2-\frac{
11911217}{59252160}a-\frac{4076357}{9875360}$  \\
 & $f_\gamma = x^4 - x^3 -4x^2+4x+1, \quad \gamma = -1.956 $ \\
 & $( -\rho^3+\rho2+\rho+1, \gamma^2- \gamma ,\rho^2-2\rho+1 )$&

\end{longtable}

\begin{longtable}[l]{l|l l }
 \caption{ramification index 3}\\
\multirow{2}{*} {1/15} 
& $ {y}^{2}+xy+y={x}^{3}+{x}^{2}-135x-660 $ \\  
& $ C=-\frac{55}{54}$  \\ 
& $\left(5 , 16 , 20\right)$\\[0.6mm]  
\hline

\multirow{2}{*} {1/10}  
& $y^2+xy+y=x^3+26-19x$ \\  
& $ C=\frac{95}{432}$\\ 
& $ \left( 10,6,15 \right)$ \\[0.6mm]  
\hline

\multirow{2}{*} {1/6i}  
&  $ y^2+xy+y=x^3+x^2-104x+101 $ \\  
&  $ C=\frac {67}{432}  $ \\    
&  $ \left( 8,7,14\right)$\\[0.6mm]  
\hline

\multirow{2}{*} {1/6ii}  
& $  y^2=x^3-x^2+16x-180 $ \\  
& $  C=1/27 $ \\  ]
& $ \left( 8,8,9\right)$  \\[0.6mm]
\hline  

\multirow{2}{*} { 5/9 }  
& $ y^2+ \left(\frac{1+\sqrt {5}}{2}\right) y = x^3+\left( \frac{1\sqrt {5}}{2}\right) {x}^{2}- \left( \frac {495 + 165\sqrt{5}}{2} \right) x-\frac{4125+1683\sqrt{5}}{2}$  \\
& $ C= -\frac{245-49\sqrt {5}}{216}$ \\   
& $\left(\sqrt{5}+3,7+3\sqrt{5},9+4\sqrt{5} \right) $ \\[0.6mm]  
\hline  

\multirow{2}{*} { 5/5 }  
& $  y^2+xy+y=x^3+x^2-110x-880$ \\  
& $  C=-\frac {35}{108} $\\  
& $ \left(5+2\sqrt {5},5+2\sqrt {5},\frac{7+3\sqrt{5}}{2} \right) $ \\[0.6mm]  
\hline

\multirow{2}{*} { 8/9 }  
& $ y^2=x^3- \left( 2\sqrt{2}+4 \right) x^2- \left( 154\sqrt {2}+231\right) x-1064\sqrt{2}-1520$ \\  
& $ C= -\frac{20+10\sqrt {2}}{27}$\\  
& $\left( 4\sqrt{2}+6,4\sqrt{2}+6,3+2\sqrt {2} \right) $\\[0.6mm]  
\hline

\multirow{2}{*}{12/3} 
& $ {y}^{2}+\sqrt {3}y={x}^{3}- \left( -\sqrt {3}+1 \right) {x}^{2}$\\  
& $  C = 0$ \\  
& $\left(4+2\sqrt {3},4+2\sqrt {3},7+4\sqrt {3} \right)$ \\[0.6mm]  
\hline

\multirow{2}{*} { 13/3i}  
& $ y^2+xy+y=x^3+x^2- \left( -495 \gamma -637 \right) x+9261\gamma+12053$ \\  
& $f_{\gamma} = x^2-x-3, \quad \gamma = 2.302$\\
& $ C=-\frac{35}{108}$\\  
& $\left(4+\sqrt{13},\frac{11+3\sqrt {13}}{2},4+\sqrt{13} \right)$ \\[0.6mm]  
\hline

\multirow{2}{*} { 13/3ii}  
& $ y^2+xy+y=x^3+x^2- \left( -495 \gamma-637 \right) x+9261\gamma+12053$ \\  
& $f_{\gamma} = x^2-x-3, \quad \gamma = -1.303  $\\
& $ C=-\frac{35}{108}$\\  
& $\left( \frac{5+\sqrt {13}}{2},22+6\sqrt {13},\frac{47+13\sqrt{13}}{2}\right)$ \\[0.6mm]  
\hline

\multirow{2}{*} { 17/36 }  
& $ y^2+xy+\gamma y=x^3-\gamma x^2-\left(19694\gamma+30770 \right)x-2145537\gamma-3350412 $ \\ 
& $ \gamma=\frac{1+\sqrt{17}}{2}$\\
& $ C= -\frac{6545+1540\sqrt{17}}{432} $ \\
& $\left(  \frac{5+1\sqrt{17}}{2}, 13+3\sqrt{17},\frac{29+7\sqrt {17}}{2} \right)$ \\[0.6mm]
\hline

\multirow{2}{*}{21/3}  
& $y^2=x^3-\left( -\frac{6 \gamma^7+7\gamma^6-30\gamma^5+4\gamma^4-108\gamma^3+90\gamma^2+42\gamma-75}{8}\right)x^2 $ \\ 
& $-\left(\frac{7\gamma^7-41 \gamma^6-20\gamma^5-126\gamma^4+180\gamma^3-108\gamma^2+51\gamma-63}{192} \right) x$ \\  
& $C= \frac{75 - 42 \gamma - 90 \gamma^2+ 36 \gamma^{3}- 4 \gamma^4+30 \gamma^5-7 \gamma^6+6\gamma^7}{144}$\\  
& $f_{\gamma} = x^8+3x^6+12x^4-9x^2+9,\quad \gamma =0.770 + 0.445 $\\ 
& $\left(\frac{5+\sqrt {7}\sqrt {3}}{2},10+2\sqrt {7}\sqrt {3},\frac{23+5\sqrt{21}}{2} \right)$ \\[0.6mm]  
\hline  

\multirow{2}{*} {28/18}  
& $ y^2+ \left( \sqrt{7}+1 \right) y=x^3- \left( \sqrt {7}+1 \right) x^2-\left( 944\sqrt{7}+2496\right) x+25532\sqrt{7}+67552 $ \\  
& $ C=\frac{1295}{576}+\frac{185}{216}\sqrt {7} $\\  
& $ \left(6+2\sqrt {7},3+\sqrt{7},8+3\sqrt {7} \right)$ \\[0.6mm]
\hline

\multirow{2}{*} {49/1} 
& $ y^2+ \left( \gamma^2+1 \right) y=x^3- \left( -\gamma^{2}-\gamma-1 \right) x^2-\left( 649\gamma^2+910\gamma+131 \right) x$ \\ 
&$-21451\gamma^{2}-21320\gamma+6760$\\
& $f_{\gamma}= x^3-x^2-2x+1,\quad \gamma=1.802   $ \\
& $C=-\frac{10 + 40\gamma + 40\gamma^2}{27}$\\  
& $\rho = 2 \cos(\frac{\pi}{7})$\\  
& $ \left( 4\rho^2+3\rho-1,4\rho^2+3\rho-1,\rho^2+\rho \right)$ \\[0.6mm]
\hline

\multirow{2}{*} { 81/1}  
&  $y^2+y=x^3-7$ \\  
&  $C=0$\\  
& $\rho = -\left( 2 \cos(\frac{5\pi}{9}) \right)^{-1} $\\
& $\left(\rho^2,\rho^2,\rho^2 \right)$ \\[0.6mm]  

\end{longtable}

\begin{longtable}[l]{l|l l }
\multicolumn{3}{l}{ramification index $e=4$}\\
\multirow{2}{*} { 8/98}  
&  $ y^2+xy+y=x^3-55146-2731x$ \\  
&  $ C=-\frac{1575}{256}$\\
&  $\left(  3+\sqrt{2},20+12\sqrt{2},21+14\sqrt {2}  \right)$ \\[0.6mm]  
\hline  

\multirow{2}{*}{8/7i}  
& $y^2=x^3- \left( 4\sqrt{2}-14 \right) x^2- \left( 32\sqrt{2}-48 \right) x$ \\  
& $C = \frac{3+6\sqrt{2}}{16}$  \\
& $ \left( 8+4\sqrt {2},4+\sqrt{2},10+6\sqrt {2}\right)$ \\[0.6mm]  
\hline  

\multirow{2}{*}{8/2i}  
& $y^2+xy=x^3- \left( 1-\sqrt{2} \right) x^2- \left( 38\sqrt{2}+51 \right) x-160\sqrt{2}-227$\\
& $ C=-\frac{87}{256}-\frac{15}{64}\sqrt{2} $ \\
& $\left(3+2\sqrt{2},7+4\sqrt {2},7+4\sqrt {2} \right)$ \\[0.6mm]
\hline  

\multirow{2}{*} {8/2ii}  
& $y^2 = 4x^3 - (1116 \sqrt{-2} +147)x - (6966\sqrt{-2} - 6859)$\\
& $ C= \frac{( -78\sqrt{-2} -123) }{2^7} $, taken from \cite{Sij2011i} \\
& $\left(3+2\sqrt{2} , 9+4\sqrt{2} ,6+4\sqrt {2} \right)$\\[0.6mm]
\hline

\multirow{2}{*}{8/7ii}  
& $y^2=x^3- \left( -4\sqrt{2}-14 \right) x^2- \left( -32\sqrt{2}-48 \right) x$ \\  
& $C = \frac{3-6\sqrt{2}}{16}$  \\
& $\left( 4+2\sqrt{2},6+2\sqrt {2},8+5\sqrt{2}\right)$ \\[0.6mm]  
\hline

\multirow{2}{*} {8/2iii}  
& $y^2+xy=x^3- \left( \sqrt{2}+1 \right) x^2- \left( -38\sqrt{2}+51 \right) x+160\sqrt{2}-227$\\ & $ C=-\frac{87}{256}+\frac{15}{64}\sqrt{2} $  \\
& $ f_{\alpha} = 5184x^4+59616x^3+171252x^2+10404x+248113, \quad \alpha= 0.1891 + 1.1341i $\\
& $\left( 5+2\sqrt{2},6+4\sqrt{2},5+2\sqrt {2} \right)$ \\[0.6mm]
\hline

\multirow{2}{*}{2624/4ii}  
& $ y^2+xy+y=x^3- \left( 1+\sqrt{2} \right) x^2- \left( -391\sqrt{2}+448 \right) x + 4342\sqrt{2}-6267$ \\ 
& $ C=-\frac{387}{256} + \frac{69}{64}\sqrt{2}$\\ 
& $ f_{\rho} = x^4-10x^3+19x^2-10x+1, \quad \rho=7.698 $ \\
& $\left(\rho,\rho, \rho+2\sqrt{\rho}+1  \right)$ \\[0.6mm]
\hline

\multirow{2}{*}{2624/4i}  
& $ y^2+xy+y=x^3- \left( 1-\sqrt{2} \right) x^2- \left( 391\sqrt{2}+448 \right) x - 4342\sqrt{2} -  6267$ \\ 
& $ C=-\frac{387}{256} - \frac{69}{64}\sqrt{2}$\\  
& $ f_{\rho_1} = x^4-10x^3+31x^2-30x+1, \quad \rho_1 = 4.965$\\
& $ f_{\rho_2} = x^5-41x^4+473x^2-1063x^2+343x-19,\quad \rho_2 = 19.181$\\
& $\left( \rho_1,\rho_2,\rho_2  \right)$ \\[0.6mm]
\hline

\multirow{2}{*} { 2304/2}  
& $  y^2-\sqrt3y=x^{3}-1$ \\  
& $  C=0$\\  
& $ \rho=3+\sqrt{2}\sqrt{3}+\sqrt{2}+\sqrt{3}$ \\
& $\left(\rho,\rho,\rho \right)$ \\[0.6mm]

\end{longtable}

\begin{longtable}[l]{l|l l }
\caption{ramification index 5}\\

\multirow{2}{*} {  5/5i }  
& $y^2+\gamma y=x^3-\gamma x^2- \left( 4217\gamma+2611 \right) x-157816\gamma-97533$ \\  
& $f_{\gamma} = x^2-x-1 , \quad \gamma= 1.6180$\\
& $C=-\frac{1083+495\sqrt{5}}{200} $ \\ 
& $f_{\gamma} = 36x^2+1962x -3299, \quad \delta=-56.132$\\
& $\left( \frac{7+1\sqrt{5}}{2},14+6\sqrt {5},16+7\sqrt {5}\right)$ \\[0.6mm]
\hline

\multirow{2}{*} {  5/180 }  
& $  y^2+xy+y=x^3-2368-334x$ \\  
& $  C=-\frac{651}{400}$ \\  
& $\left(\sqrt{5}+3,9+3\sqrt {5},\frac{21+9\sqrt {5}}{2} \right)$ \\[0.6mm]
\hline

\multirow{2}{*} {  5/5ii}  
& $y^2+\gamma y=x^3-\gamma x^2- \left( 4217\gamma+2611 \right) x-157816\gamma-97533$ \\  
& $f_{\gamma} = x^2-x-1, \quad \gamma= -0.618$\\
& $C=-\frac{1083+495\sqrt{5}}{200} $ \\  
& $\left(6+2\sqrt{5},4+\sqrt{5},\frac{17+7\sqrt{5}}{2}] \right)$ \\[0.6mm]
\hline

\multirow{2}{*} { 5/5iii }  
&  $ y^2=x^3+x^2-36x-140$ \\  
&  $ C=-\frac{6}{25}$ \\  
&  $\left(6+2\sqrt{5},6+2\sqrt{5},\frac{7+3\sqrt{5}}{2} \right)$ \\[0.6mm]
\hline

\multirow{2}{*} { 5/9 }  
& $y^2+xy+y=x^3+x^2+35x-28$\\  
& $ C=\frac{3}{100}$\\  
& $\left(\frac{9+3\sqrt{5}}{2},\frac{9+3\sqrt{5}}{2},\frac{15+5\sqrt{5}}{2} \right)$ \\[0.6mm]
\hline

\multirow{2}{*} { 725/25i}  
& $ y^2+ \left( \gamma^2+ \gamma \right) y=x^3- \left( -\gamma^3-\gamma^2+1 \right) x^2-\left( 135\gamma^3+316\gamma^{2}-136\gamma+2 \right) x $\\ 
& $-4089\gamma^3-6001\gamma^2+3228\gamma+1965 $\\  
& $f_{\gamma} = x^4-x^3-3x^2+x+1, \quad \gamma =0.738  $ \\
& $ C= \frac{-24+18\gamma-12\gamma^{2}-12\gamma^3}{25}$\\
& $f_{\rho_1} = x^4 - 3x^3 + 4x - 1, \quad \rho_1 = 2.356$\\  
& $f_{\rho_2} = x^4 - 4x^3 + 3x - 1, \quad \rho_2 = 3.811$\\  
& $\left(\rho_1 , \rho_2 ,\rho_2 \right)$ \\[0.6mm]
\hline

\multirow{2}{*} { 725/25ii}  
& $ y^2+ \left( \gamma^2+ \gamma \right) y=x^3- \left( -\gamma^3-\gamma^2+1 \right) x^2-\left( 135\gamma^3+316\gamma^{2}-136\gamma+2 \right) x $\\ 
& $-4089\gamma^3-6001\gamma^2+3228\gamma+1965 $\\  
& $f_{\gamma} = , \quad \gamma= $ \\
& $ C= \frac{-24+18\gamma-12\gamma^{2}-12\gamma^3}{25}$\\  
& $f_{\rho_1} = x^4 -  x^3 - 3x^2 + x + 1, \quad \rho_1 = 2.095$\\  
& $f_{\rho_2} = x^4 - 8x^3 + 10x^2 -x - 1, \quad \rho_2 = 6.486$\\  
& $\left(\rho_1 , \rho_2 , \rho_2 \right) $ \\[0.6mm]
\hline

\multirow{2}{*} { 1125/5 }  
& $y^2+ \left(\gamma+1 \right) y=x^3- \left( -\gamma^3+\gamma^2+1 \right)x^2- \left( 2\gamma^3-7\gamma^{2}+5\gamma+1 \right) x$\\
&$+6\gamma^3-14\gamma^2-2\gamma+12 $ \\  
&$f_{\gamma}= x^4-x^3-4x^2+4x+1 , \quad \gamma=1.338 $ \\
& $ C=0$\\  
& $f_{\rho} = x^4-3x^3-x^2+3x+1, \quad \rho = 2.956 $\\2
& $\left( \rho,\rho,\rho\right)$ \\[0.6mm]

\end{longtable}

\begin{longtable}[l]{l|l l }
\caption{ramification index 6}\\
\multirow{2}{*} {12/66i}  
& $ y^2+xy+ \left( 1-\sqrt{3} \right) y=x^3- \left( \sqrt{3}+1 \right) x^2- \left( 836-405\sqrt{3} \right) x - 4739\sqrt{3}+7704$  \\  
& $ C =\frac{53 - 387 \sqrt{3}}{54} $  \\
& $\left(3+\sqrt{3},14+6\sqrt{3},15+8\sqrt{3} \right)$\\[0.6mm]
\hline

\multirow{2}{*} {12/66ii}  
& $ y^2+xy+ \left( \sqrt{3}+1 \right) y={x}^{3}- \left( 1-\sqrt{3} \right) {x}^{2}- \left( 405\sqrt{3}+836 \right) x+4739\sqrt{3}+7704$  \\  
& $ C =\frac{53 + 387 \gamma}{54} $  \\
& $\left( 6+2\sqrt{3},5+\sqrt{3},9+4\sqrt {3}\right)$\\[0.6mm]

\end{longtable}

\begin{longtable}[l]{l|l l }
\caption{ramification index 7}\\

\multirow{2}{*} {49/91i}  
& $ y^2+xy+ay=x^3+x^2-\left( 10825\gamma^2-24436\gamma+8746 \right) x$\\
& $-995392\gamma^2+2235406\gamma-797729$ \\    
& $ C =-\frac{815}{196}+\frac{495}{49}\gamma-\frac{30}{7}\gamma^2 $ \\ 
& $ f_{\gamma} = x^3 - x^2 - 2x + 1, \quad \gamma=0.445$\\
& $ \rho = 2 \cos(\frac{\pi}{7})$\\
& $\left( \rho^2+1, 16\rho^2 + 12\rho - 8,17\rho^2+13\rho-9  \right)$\\[0.6mm]  
\hline

\multirow{2}{*} {49/91ii}  
& $ y^2+xy+ay=x^3+x^2-\left( 10825\gamma^2-24436\gamma+8746 \right) x$ \\
& $-995392\gamma^2+2235406\gamma-797729$ \\    
& $ C =-\frac{815}{196}+\frac{495}{49}\gamma-\frac{30}{7}\gamma^2 $ \\  
& $ f_{\gamma} = x^3 - x^2 - 2x + 1, \quad \gamma=-1.247$\\
& $\rho = 2\cos \left(\frac{\pi}{7} \right)$\\
& $\left( \rho^2+\rho, 5\rho^2+3\rho-2, 5\rho^2+3\rho-2 \right)$\\[0.6mm]  
\hline

\multirow{2}{*} {49/91iii}  
& $ y^2+xy+ay=x^3+x^2-\left( 10825\gamma^2-24436\gamma+8746 \right) x$ \\
& $-995392\gamma^2+2235406\gamma-797729$ \\    
& $ C = -\frac{815}{196}+\frac{495}{49}\gamma-\frac{30}{7}\gamma^2 $ \\
& $ f_{\gamma} = x^3 - x^2 - 2x + 1, \quad \gamma = 1.802$\\ 
& $\rho = 2\cos \left(\frac{\pi}{7} \right)$\\
& $\left( 2\rho^2+\rho, 2\rho^2+\rho,3\rho^2 + \rho -1 \right)$\\[0.6mm]  
\hline

\multirow{2}{*} { 49/1 }  
& $ y^2+xy+y=x^3-70-36x$\\  
& $ C = -\frac{55}{196}$\\
& $\rho = 2\cos \left(\frac{\pi}{7} \right)$\\
& $ \left( 2 \rho^2 , 2\rho^2+2\rho,4\rho^2+3\rho-2  \right)$ \\[0.6mm]

\end{longtable}

\begin{longtable}[l]{l|l l } \caption{ramification index 9} \\
\multirow{2}{*} {81/51i}  
& $ y^2= x^3- \left( -446\gamma^2-836\gamma-214 \right) x^2- \left( -375921\gamma^2-706401\gamma-199989 \right)x$ \\  
& $C=  -\frac{1309}{243}-\frac{3206}{243}\gamma-\frac{1529}{243}\gamma^{2} $\\  
& $ f_{\gamma} = x^3 - 3x -1 ,\quad  \gamma = 1.879$\\ 
& $\left(\rho^2+1, 4\rho^2 +8 \rho + 4, 5\rho^2+9\rho+3 \right)$ \\[0.6mm] 
\hline

\multirow{2}{*} {81/51ii}  
& $ y^2= x^3- \left( -446\gamma^2-836\gamma-214 \right) x^2- \left( -375921\gamma^2-706401\gamma-199989 \right)x$ \\  
& $  -\frac{1309}{243}-\frac{3206}{243}\gamma-\frac{1529}{243}\gamma^{2} $\\  
& $ f_{\gamma} = x^3 - 3x -1 ,\quad  \gamma = 1.879$\\ 
& $ \rho = 2 \cos( \frac{\pi}{9} )$\\ 
& $ \left( \rho^2 + 2 \rho + 1 , \rho^2 + 2 \rho + 1 , \rho^2 + \rho + 1     \right)$ \\[0.6mm] 
\hline

\multirow{2}{*} {81/51iii}  
& $ y^2= x^3- \left( -446\gamma^2-836\gamma-214 \right) x^2- \left( -375921\gamma^2-706401\gamma-199989 \right)x$ \\  
& $ C =  -\frac{1309}{243}-\frac{3206}{243}\gamma-\frac{1529}{243}\gamma^{2} $\\
& $ f_{\gamma} = x^3 - 3x -1 ,\quad  \gamma = -0.607 + 1.450 i  $\\ 
& $ \rho = 2 \cos( \frac{\pi}{9} )$\\ 
& $ \left( \rho^2 +2\rho+2 ,\rho^2 +2\rho+1,\rho^2 +2\rho+2    \right)$ \\[0.6mm] 
\hline

\end{longtable}
\begin{longtable}[l]{l|l l }
\caption{ramification index 11}\\
\multirow{2}{*} { 14641/1 }  
& $ y^2+y=x^3-x^2-10x-20 $ \\  
& $ C=-\frac{14}{121}$\\  
& $ \rho=2 \cos(\frac{\pi}{11})$\\  
& $ \left( \rho^3-2\rho,\rho^3-2\rho,\rho^2-1  \right)$ \\[0.6mm]

\end{longtable}

The correctness of the suggested accessory parameters is proven whenever the Fuchsian group is commensurable with a triangle group, i.e. the uniformizing differential equation is a pullback of a hypergeometric differential equation, these 25 cases can be found in \cite{Sij2011i}. The cases $1/15$ with $e=3$ can be found in Krammers article \cite{Krammer}, the case $1/6i$ was found by Elkies in \cite{Elkies}. Whenever the quaternion algebra is defined over the rational numbers, Reiter determined the accessory parameter in \cite{Reiter}. In principal it should be possible to adapt Krammers to prove at least a few more cases.

\section*{Acknowledgements}
I like to thank Stefan Reiter who convinced me that Lam\'{e} equations are interesting objects, Duco van Straten for asking and answering many questions and Jeroen Sijsling whose work builds the foundation of this article.

\bibliography{/home/hofmannj/Quellen/quellen.bib}{}

\providecommand{\bysame}{\leavevmode\hbox to3em{\hrulefill}\thinspace}
\providecommand{\MR}{\relax\ifhmode\unskip\space\fi MR }
\providecommand{\MRhref}[2]{%
  \href{http://www.ams.org/mathscinet-getitem?mr=#1}{#2}
}
\providecommand{\href}[2]{#2}
\begin{thebibliography}{EGTZ11}

\bibitem[Beu10]{Beu2010}
F.~Beukers, \emph{Recurrent sequences coming from shimura curves}, talk given
  at the Banff centre, 2010.

\bibitem[CC87]{chud87-2}
D.~V. Chudnovsky and G.~V. Chudnovsky, \emph{Computer assisted number theory
  with applications}, Number Theory (David Chudnovsky, Gregory Chudnovsky,
  Harvey Cohn, and Melvyn Nathanson, eds.), Lecture Notes in Mathematics, vol.
  1240, Springer Berlin / Heidelberg, 1987, pp.~1--68.

\bibitem[CC88]{chud88}
\bysame, \emph{Approximations and complex multiplication according to
  {R}amanujan}, Ramanujan revisited, Academic Press, Boston, MA, 1988, pp.~375
  -- 472.

\bibitem[CC89]{chud89}
\bysame, \emph{Transcendental methods and theta-functions}, Theta
  functions---{B}owdoin 1987, {P}art 2 ({B}runswick, {ME},1987), Proc. Sympos.
  Pure Math., vol.~49, Amer. Math. Soc., Providence, RI, 1989, pp.~167--232.

\bibitem[dSG]{Saint-Gervais}
H.~P. de~Saint-Gervais, \emph{Uniformisation des surfaces de riemann (french
  edition)}, ENS LSH.

\bibitem[EGTZ11]{mpc}
A.~Enge, M.~Gastineau, P.~Th\'eveny, and P.~Zimmermann, \emph{mpc --- a library
  for multiprecision complex arithmetic with exact rounding}, INRIA, 0.9 ed.,
  february 2011.

\bibitem[Elk98]{Elkies}
N.~D. Elkies, \emph{Shimura curve computations}, ANTS, 1998, pp.~1--47.

\bibitem[FBA96]{PSLQ}
H.~R.~P. Ferguson, D.~H. Bailey, and S.~Arno, \emph{Analysis of pslq, an
  integer relation finding algorithm}, Analysis \textbf{68} (1996), 351--369.

\bibitem[FK65]{FrickeKlein}
R.~Fricke and F.~Klein, \emph{Vorlesungen {\"u}ber die {T}heorie der
  automorphen {F}unktionen, volume 1}, B. G. Teubner, Stuttgart, 1965.

\bibitem[Hil08]{Hilb}
E.~Hilb, \emph{{\"U}ber {K}leinsche {T}heoreme in der {T}heorie der linearen
  {D}ifferentialgleichungen}, Mathematische Annalen \textbf{66} (1908),
  215--257.

\bibitem[Hof]{Hof1}
J.~Hofmann, \emph{A maple package for monodromy calculations, in preperation}.

\bibitem[Iha74]{Ihara}
Y.~Ihara, \emph{Schwarzian equations}, Journal of the Faculty of Science, the
  University of Tokyo. Sect. 1 A, Mathematics \textbf{21} (1974), no.~1, 97 --
  118.

\bibitem[Inc56]{ince}
E.~L. Ince, \emph{Ordinary differantial equations}, repirnt ed., Dover
  Publications, 6 1956.

\bibitem[Kat92]{Katok}
S.~Katok, \emph{{F}uchsian {G}roups ({C}hicago {L}ectures in {M}athematics)}, 1
  ed., University Of Chicago Press, 8 1992.

\bibitem[Kra96]{Krammer}
D.~Krammer, \emph{An example of an arithmetic {F}uchsian group}, J. reine
  angew. Math. \textbf{473} (1996), 69--85.

\bibitem[Poo07]{Poole}
E.~G.~C. Poole, \emph{Introduction to the theory of linear differential
  equations}, Pierides Press, 3 2007.

\bibitem[Rei09]{Reiter}
S.~Reiter, \emph{Halphen's transform and middle convolution}, ArXiv e-prints
  (2009).

\bibitem[Sija]{Sij2011i}
J.~Sijsling, \emph{Arithmetic (1;e)-curves and {B}elyi maps}, Math.Comp., to
  appear.

\bibitem[Sijb]{Sij2011}
J.~Sijsling, \emph{Canonical models of arithmetic (1;e)-curves}.

\bibitem[Sijc]{Sij-tools}
J.~Sijsling, \emph{Magma programs for arithmetic pointed tori}.

\bibitem[Tak83]{Tak83}
K.~Takeuchi, \emph{Arithmetic {F}uchsian groups with signature (1;e)}, J. Math.
  Soc. Japan \textbf{35} (1983), 381--407.

\bibitem[Vig80]{vign1980}
M.~F. Vign{\'e}ras, \emph{Arithm{\'e}tique des alg{\`e}bres de quaternions.},
  Lecture notes in mathematics, no. 800, Springer-Verlag, 1980.

\bibitem[Yos87]{Yoshida1987}
M.~Yoshida, \emph{Fuchsian differential equations}, Friedrich Vieweg und Son,
  1987.

\end{thebibliography}
\bibliographystyle{amsalpha}

\end{document}